\documentclass[11pt]{amsart}
\usepackage{amsmath,amssymb,amsthm,mathtools}
\usepackage{enumitem}
\usepackage{booktabs}
\usepackage{caption}
\usepackage{pifont}
\usepackage{xcolor}
\DeclareRobustCommand{\commonaddress}{3, Bunkyo-cho, Hirosaki-shi, Aomori, 036-8561, JAPAN}
\usepackage{url}
\usepackage[colorlinks=true,linkcolor=blue,citecolor=blue,urlcolor=blue]{hyperref}

\mathtoolsset{showonlyrefs=true}
\allowdisplaybreaks

\newtheorem{theorem}{Theorem}[section]
\newtheorem{proposition}[theorem]{Proposition}
\newtheorem{lemma}[theorem]{Lemma}
\newtheorem{corollary}[theorem]{Corollary}
\newtheorem{definition}[theorem]{Definition}

\newtheorem{fact}[theorem]{Fact}

\newcommand{\one}{\mathbf 1}
\newcommand{\wt}{\operatorname{wt}}

\title[Maximal Triply Even Binary Codes]
{Maximal Triply Even Binary Codes from Elliptic Orthogonal Graphs}
\author[Kasai and Nogata]{{Yusuke Kasai$^\dagger$ and Yuto Nogata$^*$\\
\normalfont\small \commonaddress}}
\thanks{{$^*$ Graduate School of Science and Technology,
Hirosaki University (corresponding author, h26ds201@hirosaki-u.ac.jp).}}
\thanks{{$^\dagger$ Graduate School of Science and Technology,
Hirosaki University. Email: h25ms106@hirosaki-u.ac.jp}}

\subjclass[2020]{{Primary 94B05; Secondary 05E30, 51A50}}
\keywords{{Maximal triply even binary code, orthogonal polar space,
strongly regular graph}}
\date{}

\begin{document}

\begin{abstract}
For a finite set \(X\), a binary linear code
\(C\subseteq\mathbb F_2^X\) is called a triply even binary code if
\(\wt(c)\equiv0\pmod8\) for every \(c\in C\).
Let \(q\) be an odd prime power and \(n\) be an integer with
\(n\ge2\). Let \(E^-_{2n+2}(q)\) denote the binary code spanned by all
pairwise sums of rows of the adjacency matrix of the elliptic orthogonal
graph arising from a \((2n+2)\)-dimensional elliptic quadratic space over
\(\mathbb F_q\). We prove that
\(E^-_{2n+2}(q)\) is a triply even binary code if and only if
\(n=2\) and \(q\equiv3\pmod4\). In this case, we also prove that
\(E^-_6(q)\) is a maximal triply even binary code.
\end{abstract}

\maketitle

\section{Introduction}

Let $X$ be a finite set. A subspace
$C\subseteq\mathbb F_2^X$ is called a
\emph{binary linear code}.
If $|X|=N$ and $\dim C=K$, then $C$ is called an
\emph{$[N,K]$ binary linear code}.
For $c=(c_x)_{x\in X}\in C$, define the
\emph{Hamming weight} by
$
 \wt(c):=\#\{x\in X\mid c_x=1\}.
$
We begin by recalling the definition of a maximal triply even binary code.

\begin{definition}

A binary linear code $C\subseteq\mathbb F_2^X$ is called a
\emph{triply even binary code} if
$\wt(c)\equiv0\pmod8$ for every $c\in C$.
A triply even binary code $C$ is called a
\emph{maximal triply even binary code} if there is no triply even
binary code $C'$ satisfying $C\subsetneq C'\subsetneq\mathbb F_2^X$.
\end{definition}

We recall the following known result.

\begin{fact}[{\cite[Section~4]{BetsumiyaMunemasa}}]\label{fac:betsumiya_2TEC}
 Let $m\ge4$. The code generated by the rows of the adjacency matrix
of the triangular graph $T_m$ is a maximal triply even binary code if and only if
$m\equiv2\pmod4$.
\end{fact}
Throughout this paper, let $q$ be an odd prime power {and $n$ be an
integer with $n\ge2$}.
Motivated by Fact~\ref{fac:betsumiya_2TEC}, we consider the binary code
$E^-_{2n+2}(q)$ constructed from the elliptic orthogonal graph
$\Gamma(O^-_{2n+2}(q))$ using pairwise sums of rows of its adjacency matrix.
The graph and the code are defined in \S~2.1.

Our main result is the following.

\begin{theorem}\label{thm:main}
Let $q$ be an odd prime power and $n$ be an integer with $n\ge2$.
Then $E^-_{2n+2}(q)$ is a triply even binary code if and only if
$n=2$ and $q\equiv3\pmod4$. Moreover, if $q\equiv3\pmod4$, then
$E^-_6(q)$ is a
$
  \bigl[(q+1)(q^3+1),\ q(q^2-q+1)\bigr]
$
maximal triply even binary code.
\end{theorem}

\begin{proof}
The first assertion follows from Proposition~\ref{prop:o6-triply-even} and
Corollary~\ref{cor:higher-rank-not-triply-even}. The second assertion follows from
Corollary~\ref{cor:o6-maximal}.
\end{proof}

This paper is organized as follows:
\begin{itemize}
  \item \S~2 introduces the notation and lemmas used in this paper.
  \item \S~3 proves that \(E^-_{2n+2}(q)\) is a triply even binary code
  if and only if \(n=2\) and \(q\equiv3\pmod4\).
  \item \S~4 proves that \(E^-_6(q)\) is a maximal triply even binary code
  when \(q\equiv3\pmod4\).
\end{itemize}

\section{Notation and preliminaries}

 {In this section, we introduce the notation and preliminary
results used in this paper.}
{We also outline the strategy for \S~3 and \S~4.}
\begin{itemize}
 \item {In \S~2.1, we define the elliptic orthogonal polar space
 $O^-_{2n+2}(q)$, the graph $\Gamma(O^-_{2n+2}(q))$, and the code
 $E^-_{2n+2}(q)$ considered in this paper. We also recall some properties of
 $O^-_6(q)$.}

 \item {In \S~2.2, we show that $E^-_{2n+2}(q)$ is a triply even
 binary code if and only if \textup{(E1)}-\textup{(E3)} hold. We also give a
 sufficient condition for its maximality.}

 \item {In \S~2.3, we prepare lemmas on strongly regular graphs
 for verifying \textup{(E1)}-\textup{(E3)}.}
\end{itemize}

Let $\Gamma$ be a finite simple graph with vertex set $X$. We write
$x\sim y$ when $x$ and $y$ are adjacent, and set
$
 N(x):=\{y\in X\mid x\sim y\}.
$

\subsection{\texorpdfstring{{Elliptic orthogonal polar space, graph, and code}}{Elliptic orthogonal polar space, graph, and code}}

{In this subsection, we define $O^-_{2n+2}(q)$,
$\Gamma(O^-_{2n+2}(q))$, and $E^-_{2n+2}(q)$.}

\begin{definition}
Let \(\mathcal P\) be a set whose elements are called \emph{points}, and
\(\mathcal L\) be a collection of subsets of \(\mathcal P\) whose members are
called \emph{lines}. The pair $(\mathcal P,\mathcal L)$ is called a
\emph{partial linear space} if the following conditions hold.
\begin{enumerate}[label=\textup{(P\arabic*)}]
 \item Any two distinct points lie on at most one line.
 \item Every line contains at least two points.
\end{enumerate}
Two distinct points are called \emph{collinear} if they lie on a common line.
\end{definition}

\begin{definition}

A partial linear space $(\mathcal P,\mathcal L)$ is called a
\emph{generalized quadrangle} if it satisfies the following condition.
\begin{enumerate}[label=\textup{(GQ)}]
 \item For every $P\in\mathcal P$ and $\ell\in\mathcal L$ with
 $P\notin\ell$, there exists a unique point on $\ell$ collinear with $P$.
\end{enumerate}
If every line contains $s+1$ points and every point lies on $t+1$ lines,
then $(\mathcal P,\mathcal L)$ is said to have \emph{order} $(s,t)$.
\end{definition}

\begin{definition}
Let $(\mathcal P,\mathcal L)$ be a generalized quadrangle.
Three distinct points $P_1,P_2,P_3\in\mathcal P$ form a
\emph{triad} if they are pairwise non-collinear.
A point $Z\in\mathcal P$ is called a \emph{center} of the triad if
$Z$ is collinear with each of $P_1,P_2,$ and $P_3$.
\end{definition}

\begin{fact}[{\cite[Theorem~5.11]{Pech}}]
\label{fac:gq-triad-centres}
Let $(\mathcal P,\mathcal L)$ be a generalized quadrangle of order
$(s,t)$. Then $s^2\ge t$. Moreover, equality holds if and only if every
triad in $(\mathcal P,\mathcal L)$ has the same number $s+1$ of centers.
\end{fact}

A \emph{quadratic form} \(f\) on a vector space \(V\)
over \(\mathbb F_q\) is a mapping \(f\colon V\to\mathbb F_q\) satisfying
\begin{equation}\label{eq:quadratic-polar}
 f(\alpha v)=\alpha^2f(v),\quad
 f(u+v)=f(u)+f(v)+B_f(u,v)
\end{equation}
for any \(u,v\in V\) and \(\alpha\in\mathbb F_q\), where \(B_f\) is a symmetric
bilinear form.
We call $(V,f)$ a \emph{quadratic space}.

The \emph{radical} of $f$ is defined by
$
 \operatorname{Rad}f
 :=\{x\in V\mid f(x)=0,\ B_f(x,y)=0\text{ for all }y\in V\}
$.
The quadratic space $(V,f)$ is called
\emph{nondegenerate} if
$\operatorname{Rad}f=\{0\}$.
A subspace $W\subseteq V$ is called
\emph{singular} if $f(w)=0$ for every $w\in W$.
If $\dim_{\mathbb F_q} W=r$, then $W$ is also called a
\emph{singular $r$-space}.

Let \((V,f)\) be a \(2m\)-dimensional nondegenerate quadratic space.
Every maximal singular subspace of $(V,f)$
has dimension either $m-1$ or $m$.
If \((V,f)\) has maximal singular subspaces of dimension \(m-1\), then
it is called an \emph{elliptic quadratic space}.

Let $V$ be a nonzero finite-dimensional vector space over $\mathbb F_q$.
The projective space associated with $V$ is
$\operatorname{PG}(V)
:=\{\langle x\rangle_{\mathbb F_q}\mid x\in V\setminus\{0\}\}$.
Its elements are called \emph{points}.
For each two-dimensional subspace $W$ of $V$, the set
$\ell_W:=\operatorname{PG}(W)
=\{\langle w\rangle_{\mathbb F_q}\mid w\in W\setminus\{0\}\}$
is called a \emph{projective line}.
We regard a projective line as a set of points and write
$P\in\ell_W$ when the point $P$ lies on $\ell_W$.

\begin{definition}
Let $q$ be an odd prime power and $n$ be an integer with $n\ge2$.
Let $(V,f)$ be a $(2n+2)$-dimensional elliptic quadratic space
over $\mathbb F_q$. Define
\begin{itemize}
 \item $\mathcal P_{2n+1}
 :=\{\langle x\rangle_{\mathbb F_q}\mid
 x\in V\setminus\{0\},\ f(x)=0\}$,
 \item $\mathcal L_{2n+1}
 :=\{\ell_W\mid W\subsetneq V
 \text{ is a singular $2$-subspace}\}$.
\end{itemize}
The partial linear space $(\mathcal P_{2n+1},\mathcal L_{2n+1})$
is called the \emph{elliptic orthogonal polar space}
associated with $(V,f)$ and is denoted by $O^-_{2n+2}(q)$.
\end{definition}

\begin{definition}
 We define the \emph{elliptic orthogonal graph} $\Gamma(O^-_{2n+2}(q))$
as the graph with vertex set $\mathcal P_{2n+1}$, where two distinct vertices
are adjacent if and only if they are collinear.
\end{definition}

\begin{definition}
{For $P\in\mathcal P_{2n+1}$, let
$r_P\in\mathbb F_2^{\mathcal P_{2n+1}}$
denote the row corresponding to $P$ in the adjacency matrix of
$\Gamma(O^-_{2n+2}(q))$.}
{Fix $P_0\in\mathcal P_{2n+1}$. For
$P\in\mathcal P_{2n+1}\setminus\{P_0\}$, put $g_P:=r_P+r_{P_0}$.}
We define
\begin{equation}\label{eq:elliptic-code}
 E^-_{2n+2}(q)
 :=\langle r_P+r_Q\mid P,Q\in\mathcal P_{2n+1}\rangle_{\mathbb F_2}
 =\langle g_P\mid
 P\in\mathcal P_{2n+1}\setminus\{P_0\}\rangle_{\mathbb F_2}.
\end{equation}

\end{definition}

{The space $O^-_6(q)$ has useful properties for studying
$E^-_6(q)$.}

\begin{lemma}[{\cite[Theorem~2.3.5 and Table~2.1]
{BrouwerVanMaldeghem}}]
\label{lem:o6-generalized-quadrangle}

Let $q$ be an odd prime power. Then the elliptic orthogonal polar space
$O^-_6(q)$ is a generalized quadrangle of order $(q,q^2)$. Hence, the following properties hold.
\begin{enumerate}
 \item[\textup{(P1)$'$}]
 If distinct $P_1,P_2\in\mathcal P_5$ satisfy
 $P_1,P_2\in\ell\cap m$ for $\ell,m\in\mathcal L_5$, then $\ell=m$.

 \item[\textup{(P2)$'$}]
 Every $\ell\in\mathcal L_5$ contains exactly $q+1$ points.

 \item[\textup{(GQ)$'$}]
 If $P\in\mathcal P_5$ and $\ell\in\mathcal L_5$ satisfy $P\notin\ell$,
 then $|N(P)\cap\ell|=1$.
\end{enumerate}
\end{lemma}

\begin{lemma}\label{lem:gq-triangles}

Let $P_1,P_2\in\mathcal P_5$ be distinct points such that $P_1\sim P_2$,
and $\ell_{12}$ be the unique line containing $P_1$ and $P_2$.
Then $N(P_1)\cap N(P_2)=\ell_{12}\setminus\{P_1,P_2\}$.
\end{lemma}

\begin{proof}
Let $P_3\in N(P_1)\cap N(P_2)$. Then $P_1\sim P_3$ and
$P_2\sim P_3$, and since the graph is simple,
$P_3\notin\{P_1,P_2\}$. Suppose that $P_3\notin\ell_{12}$.
Since $P_1,P_2\in\ell_{12}$ and both are adjacent to $P_3$, we have
$P_1,P_2\in N(P_3)\cap\ell_{12}$, contradicting \textup{(GQ)$'$}.
Hence $P_3\in\ell_{12}\setminus\{P_1,P_2\}$.

Conversely, if $P_3\in\ell_{12}\setminus\{P_1,P_2\}$, then
$P_1\sim P_3$ and $P_2\sim P_3$ by the definition of adjacency.
Thus $P_3\in N(P_1)\cap N(P_2)$.
\end{proof}

\subsection{Triply even binary codes and maximality}

{In this subsection, we prepare the conditions used later
to determine when $E^-_{2n+2}(q)$ is a triply even binary code and when it is maximal.}

For \(c,d\in\mathbb F_2^X\), put
\(c*d:=(c_xd_x)_{x\in X}\) and
\(\langle c,d\rangle:=\sum_{x\in X}c_xd_x\), where
\(\langle\cdot,\cdot\rangle\) denotes the standard inner product on
\(\mathbb F_2^X\). For \(C,D\subseteq\mathbb F_2^X\), define
\(C^\perp:=\{u\in\mathbb F_2^X\mid
\langle u,c\rangle=0\text{ for every }c\in C\}\) and
\(C*D:=\langle c*d\mid c\in C,\ d\in D\rangle_{\mathbb F_2}\).

\begin{lemma}[{\cite[Lemma~1]{BetsumiyaMunemasa}}]
\label{lem:generator-criterion}

Let
$C=\langle g_1,\ldots,g_d\rangle_{\mathbb F_2}
\subseteq\mathbb F_2^X$.
Then $C$ is a triply even binary code if and only if the following three
conditions hold.
\begin{enumerate}[label=\textup{(G\arabic*)}]
 \item $\wt(g_i)\equiv0\pmod8$ for every $i$.
 \item $\wt(g_i*g_j)\equiv0\pmod4$ for every $i<j$.
 \item $\wt(g_i*g_j*g_k)\equiv0\pmod2$ for every $i<j<k$.
\end{enumerate}
\end{lemma}

\begin{corollary}\label{cor:elliptic-generator-criterion}
 $E^-_{2n+2}(q)$ is a triply even binary code if and only if the
following conditions hold.
\begin{itemize}
 \item[\textup{(E1)}]
 $\wt(g_P)\equiv0\pmod8$ for every
 $P\in\mathcal P_{2n+1}\setminus\{P_0\}$.

 \item[\textup{(E2)}]
 $\wt(g_P*g_R)\equiv0\pmod4$ for all distinct
 $P,R\in\mathcal P_{2n+1}\setminus\{P_0\}$.

 \item[\textup{(E3)}]
 $\wt(g_P*g_R*g_T)\equiv0\pmod2$ for all distinct
 $P,R,T\in\mathcal P_{2n+1}\setminus\{P_0\}$.
\end{itemize}
\end{corollary}

{In \S~3, we determine when conditions
\textup{(E1)}-\textup{(E3)} hold.}

\begin{proposition}[{\cite[Lemmas 3, 8]{BetsumiyaMunemasa}}]\label{prop:triply-even-square-maximal}
 Let $C\subseteq\mathbb F_2^X$ be a triply even binary code. Then the following
statements hold.
\begin{enumerate}[label=\textup{(M\arabic*)}]
 \item $C*C\subseteq C^\perp$.
 \item If $(C*C)^\perp=C$, then $C$ is a maximal triply even binary code.
\end{enumerate}

\end{proposition}

\begin{corollary}\label{cor:elliptic-maximality-criterion}
{
Suppose that $E^-_{2n+2}(q)$ is a triply even binary code. If
\begin{equation}\label{eq:elliptic-maximality-condition}
 E^-_{2n+2}(q)*E^-_{2n+2}(q)
 =E^-_{2n+2}(q)^\perp,
\end{equation}
then $E^-_{2n+2}(q)$ is a maximal triply even binary code.
}
\end{corollary}

{In \S~4, we prove that
\eqref{eq:elliptic-maximality-condition} holds for $n=2$.}

\subsection{\texorpdfstring{{Strongly regular graphs}}{Strongly regular graphs}}
{In this subsection, we prepare lemmas on strongly regular graphs
for verifying \textup{(E1)}-\textup{(E3)}.}

\begin{definition}\label{def:strongly-regular-graph}

Let $v,k,\lambda,$ and $\mu$ be integers such that
$0\le k\le v-1$ and $\lambda,\mu\ge0$.
Then $\Gamma$ is called a \emph{strongly regular graph} with parameters
$(v,k,\lambda,\mu)$ if $|X|=v$, $\Gamma$ is $k$-regular, and for every
pair of distinct vertices $x,y\in X$,
\[
 |N(x)\cap N(y)|=
 \begin{cases}
  \lambda,&x\sim y,\\
  \mu,&x\not\sim y.
 \end{cases}\]
\end{definition}

We use the notation $\operatorname{srg}(v,k,\lambda,\mu)$ for a strongly
regular graph with parameters $(v,k,\lambda,\mu)$.

\begin{lemma}\label{lem:srg-row-difference-weight}

Let $\Gamma$ be an $\operatorname{srg}(v,k,\lambda,\mu)$, and $r_u\in \mathbb{F}_{2}^{X}$ denote the row of its adjacency matrix corresponding to $u\in X$.
For any distinct vertices $x,y\in X$, we have
\[
 \wt(r_x+r_y)=
 \begin{cases}
  2(k-\lambda),&x\sim y,\\
  2(k-\mu),&x\not\sim y.
 \end{cases}
\]
\end{lemma}

\begin{proof}
By the definition of a strongly regular graph, we have
\[
 \wt(r_x+r_y)
 =|N(x)|+|N(y)|-2|N(x)\cap N(y)|
 =
 \begin{cases}
  2(k-\lambda),&x\sim y,\\
  2(k-\mu),&x\not\sim y.
 \end{cases}
\]
\end{proof}

\begin{proposition}\label{prop:srg-generator-criterion-s1-s2}
Let $0< k< v-1$, $\Gamma$ be an $\operatorname{srg}(v,k,\lambda,\mu)$, and fix $w\in X$.
Then the following two conditions are equivalent.
\begin{enumerate}[label=\textup{(S\arabic*)}]
 \item $\wt(r_x+r_w)\equiv0\pmod8$ for every
 $x\in X\setminus\{w\}$, and
 $
  \wt((r_x+r_w)*(r_y+r_w))\equiv0\pmod4
 $
 for every pair of distinct vertices $x,y\in X\setminus\{w\}$.
 \item $4\mid(k-\lambda)$ and $4\mid(k-\mu)$.
\end{enumerate}
\end{proposition}

\begin{proof}
Since $0<k<v-1$, the vertex $w$ has both a neighbor and a nonneighbor.
Assume \textup{(S1)}. By Lemma~\ref{lem:srg-row-difference-weight},
both $2(k-\lambda)$ and $2(k-\mu)$ are divisible by $8$. Hence
$4\mid(k-\lambda)$ and $4\mid(k-\mu)$.

Conversely, assume \textup{(S2)}.
By Lemma~\ref{lem:srg-row-difference-weight},
$8\mid\wt(r_a+r_b)$ for all distinct $a,b\in X$.
For distinct $x,y\in X\setminus\{w\}$,
\[
 2\wt((r_x+r_w)*(r_y+r_w))
 =\wt(r_x+r_w)+\wt(r_y+r_w)-\wt(r_x+r_y).
\]
Hence
$4\mid\wt((r_x+r_w)*(r_y+r_w))$.
\end{proof}

\begin{fact}[{\cite[Section~2.6.3]{BrouwerVanMaldeghem}}]
\label{fac:elliptic-srg-parameters}

Let $q$ be an odd prime power and $n$ be an integer with $n\ge2$. Then
$\Gamma(O^-_{2n+2}(q))$ is an $\operatorname{srg}(v,k,\lambda,\mu)$, where
\[
\begin{aligned}
 v&=\frac{(q^n-1)(q^{n+1}+1)}{q-1},&
 k&=q\frac{(q^{n-1}-1)(q^n+1)}{q-1},\\
 \lambda&=q^2\frac{(q^{n-2}-1)(q^{n-1}+1)}{q-1}+q-1,&
 \mu&=\frac{(q^{n-1}-1)(q^n+1)}{q-1}.
\end{aligned}
\]
The parameters $k,\lambda,$ and $\mu$ are all even.
In particular, when $n=2$, the parameters of $\Gamma(O^-_6(q))$ are
\begin{equation}\label{eq:o6-srg-parameters}
 (v,k,\lambda,\mu)
 =\bigl((q+1)(q^3+1),\ q(q^2+1),\ q-1,\ q^2+1\bigr).
\end{equation}
\end{fact}

For distinct $a,b,c\in X$, put
$
 \tau(a,b,c):=|N(a)\cap N(b)\cap N(c)|.
$

\begin{lemma}\label{lem:triple-intersection-transfer}
 Let $\Gamma$ be an $\operatorname{srg}(v,k,\lambda,\mu)$, and suppose that
$k$, $\lambda$, and $\mu$ are all even.
For any four distinct vertices $x,y,z,w\in X$, we have
{
\begin{equation}\label{eq:triple-intersection-transfer}
\begin{aligned}
&\wt((r_x+r_w)*(r_y+r_w)*(r_z+r_w))\\
&\equiv \tau(x,y,z)+\tau(w,y,z)+\tau(w,x,z)+\tau(w,x,y)\pmod2.
\end{aligned}
\end{equation}
}
\end{lemma}

\begin{proof}
Since $r_w*r_w=r_w$,
{
\[
\begin{aligned}
&(r_x+r_w)*(r_y+r_w)*(r_z+r_w)\\
&=r_x*r_y*r_z+r_w*r_y*r_z+r_x*r_w*r_z+r_x*r_y*r_w\\
&\quad+r_w*r_x+r_w*r_y+r_w*r_z+r_w.
\end{aligned}
\]
}
For $u_1,\ldots,u_t\in \mathbb F_2^X$,
$
 \wt\left(\sum_{i=1}^t u_i\right)
 \equiv\sum_{i=1}^t\wt(u_i)\pmod2
$. If $a,b\in X$ are distinct, then
$\wt(r_a*r_b)=|N(a)\cap N(b)|\in\{\lambda,\mu\}$, and
$\wt(r_w)=k$. Hence
{
\[
\begin{aligned}
 &\wt((r_x+r_w)*(r_y+r_w)*(r_z+r_w))\\
 & \equiv  \tau(x,y,z)+\tau(w,y,z)+\tau(w,x,z)+\tau(w,x,y)\\
&\quad+\wt(r_w*r_x)+\wt(r_w*r_y)+\wt(r_w*r_z)+\wt(r_w)\\
 &\equiv \tau(x,y,z)+\tau(w,y,z)+\tau(w,x,z)+\tau(w,x,y)\pmod2.
\end{aligned}
\]
}
\end{proof}

\section{\texorpdfstring{{When $E^-_{2n+2}(q)$ is a triply even binary code}}{When E^-_{2n+2}(q) is a triply even binary code}}

In this section, we prove that $E^-_{2n+2}(q)$ is a triply even binary code if and only if
$n=2$ and $q\equiv3\pmod4$.
By Corollary~\ref{cor:elliptic-generator-criterion}, it suffices to
determine when conditions \textup{(E1)}, \textup{(E2)}, and \textup{(E3)}
hold. We first determine when \textup{(E1)} and \textup{(E2)} hold.

\begin{lemma}\label{lem:elliptic-generator-criterion-e1-e2}
\textup{(E1)} and \textup{(E2)} hold if and only if
$q\equiv3\pmod4$.
\end{lemma}

\begin{proof}
By Fact~\ref{fac:elliptic-srg-parameters},
$
 k-\lambda=q^{2n-1}+1$, $k-\mu=(q^{n-1}-1)(q^n+1)$.
Applying Proposition~\ref{prop:srg-generator-criterion-s1-s2} with
$w=P_0$, conditions \textup{(E1)} and \textup{(E2)} hold if and only if
\begin{equation}\label{eq:elliptic-divisibility-i-ii}
 4\mid(k-\lambda),
 \qquad
 4\mid(k-\mu).
\end{equation}

Since $q^2\equiv1\pmod4$, we have $q^{2n-1}\equiv q\pmod4$. Hence
\begin{equation}\label{eq:k-lambda-divisibility}
 4\mid(k-\lambda)\iff q\equiv3\pmod4.
\end{equation}
Thus \eqref{eq:elliptic-divisibility-i-ii} implies
$q\equiv3\pmod4$.

Conversely, assume $q\equiv3\pmod4$. If $n$ is even, then
$q^{n-1}-1\equiv2\pmod4$ and $q^n+1\equiv2\pmod4$. If $n$ is odd, then
$q^{n-1}-1\equiv0\pmod4$. Hence $4\mid(k-\mu)$ in either case, and
\eqref{eq:k-lambda-divisibility} gives $4\mid(k-\lambda)$.
Therefore \eqref{eq:elliptic-divisibility-i-ii} holds.
\end{proof}

It remains to determine when \textup{(E3)} holds.
In \S~3.1, we prove that \textup{(E3)} holds when $n=2$, whereas in
\S~3.2, we prove that it does not hold when $n\ge3$.

\subsection{\texorpdfstring{The case $n=2$}{The case n=2}}
In this subsection, we prove that \textup{(E3)} holds for $E^-_6(q)$.

\begin{corollary}\label{cor:o6-triad-centres}
Every triad in $O^-_6(q)$ has exactly $q+1$ centers.
\end{corollary}

\begin{proof}
By Lemma~\ref{lem:o6-generalized-quadrangle},
$O^-_6(q)$ is a generalized quadrangle of order $(q,q^2)$.
By Fact~\ref{fac:gq-triad-centres}, every triad has exactly
$q+1$ centers.
\end{proof}

For three distinct points $P_1,P_2,P_3\in\mathcal P_5$, put
\[
 e(P_1,P_2,P_3)
 :=
 \#\left\{
 \{i,j\}\in\binom{\{1,2,3\}}2
 \ \middle|\ P_i\sim P_j
 \right\}.
\]

\begin{proposition}\label{prop:o6-triple-neighbours}
For three distinct points $P_1,P_2,P_3\in\mathcal P_5$,
$\tau(P_1,P_2,P_3)$ is given by the following table.
\begin{center}
\captionof{table}{Values of $e(P_1,P_2,P_3)$ and
$\tau(P_1,P_2,P_3)$}
\label{tab:o6-triple-neighbours}
\begin{tabular}{c@{\qquad}c}
\toprule
$e(P_1,P_2,P_3)$ & $\tau(P_1,P_2,P_3)$\\
\midrule
$0$ & $q+1$\\
$1$ & $1$\\
$2$ & $0$\\
$3$ & $q-2$\\
\bottomrule
\end{tabular}
\end{center}
In particular,
$\tau(P_1,P_2,P_3)\equiv e(P_1,P_2,P_3)\pmod2$.
\end{proposition}

\begin{proof}
Assume $e(P_1,P_2,P_3)=0$. Then $\{P_1,P_2,P_3\}$ is a triad in $O^-_6(q)$.
By Corollary~\ref{cor:o6-triad-centres}, it has exactly $q+1$ centers,
and hence
$
 \tau(P_1,P_2,P_3)=q+1$.

Assume $e(P_1,P_2,P_3)=1$. We may take
$
 P_1\sim P_2$, $P_1\not\sim P_3$, and $P_2\not\sim P_3$.
Let $\ell$ be the line containing $P_1$ and $P_2$.
Then $P_3\notin\ell$. By \textup{(GQ)$'$}, there exists a unique point
$Y\in\ell$ collinear with $P_3$.
Since $P_1,P_2\notin N(P_3)$, we have
$Y\notin\{P_1,P_2\}$, and hence
$Y\in N(P_1)\cap N(P_2)\cap N(P_3)$.
{
If $Y'\in N(P_1)\cap N(P_2)\cap N(P_3)$, then
Lemma~\ref{lem:gq-triangles} gives $Y'\in\ell$.
Since $Y'\sim P_3$, the uniqueness of $Y$ implies $Y'=Y$.
}
Hence
$
 N(P_1)\cap N(P_2)\cap N(P_3)=\{Y\},
$
and therefore $\tau(P_1,P_2,P_3)=1$.

Assume $e(P_1,P_2,P_3)=2$. We may take
$
 P_1\sim P_2$,
 $P_2\sim P_3$,
 $P_1\not\sim P_3$.
Let $\ell_{12}$ and $\ell_{23}$ be the lines containing
$P_1,P_2$ and $P_2,P_3$, respectively.
Since $P_1\not\sim P_3$, these lines are distinct, and \textup{(P1)$'$} gives
$
 \ell_{12}\cap\ell_{23}=\{P_2\}$.
If $Y\in N(P_1)\cap N(P_2)\cap N(P_3)$, then
Lemma~\ref{lem:gq-triangles} gives
$Y\in\ell_{12}\cap\ell_{23}$, and hence $Y=P_2$.
This is impossible because $P_2\notin N(P_2)$.
Therefore $\tau(P_1,P_2,P_3)=0$.

Assume $e(P_1,P_2,P_3)=3$.
Let $\ell$ be the line containing $P_1$ and $P_2$.
Since $P_3\in N(P_1)\cap N(P_2)$,
Lemma~\ref{lem:gq-triangles} gives $P_3\in\ell$.
Hence $P_1,P_2,$ and $P_3$ lie on the same line $\ell$. Let $Y\in N(P_1)\cap N(P_2)\cap N(P_3)$.
Since $Y\sim P_1$ and $Y\sim P_2$,
Lemma~\ref{lem:gq-triangles} gives $Y\in\ell$.
Thus
$
 N(P_1)\cap N(P_2)\cap N(P_3)
 \subseteq \ell\setminus\{P_1,P_2,P_3\}.
$
Conversely, every point of
$\ell\setminus\{P_1,P_2,P_3\}$ is adjacent to each of
$P_1,P_2,$ and $P_3$. Hence
$
 N(P_1)\cap N(P_2)\cap N(P_3)
 =\ell\setminus\{P_1,P_2,P_3\}.
$
By \textup{(P2)$'$},
$|\ell|=q+1$, and therefore
$\tau(P_1,P_2,P_3)=q-2$.

Since $q$ is odd,
$
 \tau(P_1,P_2,P_3)\equiv e(P_1,P_2,P_3)\pmod2$.
\end{proof}

\begin{lemma}\label{lem:o6-generator-criterion-iii}
\textup{(E3)} holds for $E^-_6(q)$.
\end{lemma}

\begin{proof}
Let $P_1,P_2,P_3\in\mathcal P_5\setminus\{P_0\}$ be distinct.
For distinct $i,j\in\{0,1,2,3\}$, put
$\varepsilon_{i,j}=1$ if $P_i\sim P_j$, and
$\varepsilon_{i,j}=0$ otherwise.
Then, for distinct $i,j,k\in\{0,1,2,3\}$,
$
e(P_i,P_j,P_k)=\varepsilon_{i,j}+\varepsilon_{i,k}+\varepsilon_{j,k}$.

By Fact~\ref{fac:elliptic-srg-parameters}, the parameters
$k,\lambda,$ and $\mu$ of $\Gamma(O^-_6(q))$ are all even.
Hence, by Lemma~\ref{lem:triple-intersection-transfer} and
Proposition~\ref{prop:o6-triple-neighbours},
{
\[
\begin{aligned}
 &\wt(g_{P_1}*g_{P_2}*g_{P_3})\\
 &\equiv \tau(P_1,P_2,P_3)+\tau(P_0,P_2,P_3)+\tau(P_0,P_1,P_3)+\tau(P_0,P_1,P_2)\\
 &\equiv e(P_1,P_2,P_3)+e(P_0,P_2,P_3)+e(P_0,P_1,P_3)+e(P_0,P_1,P_2)\\
 &=2\sum_{0\le i<j\le3}\varepsilon_{i,j}
 \equiv0\pmod2.
\end{aligned}
\]
}
Thus \textup{(E3)} holds.
\end{proof}

\begin{proposition}\label{prop:o6-triply-even}
$E^-_6(q)$ is a triply even binary code if and only if
$q\equiv3\pmod4$.
\end{proposition}

\begin{proof}
{The assertion follows from
Corollary~\ref{cor:elliptic-generator-criterion},
Lemmas~\ref{lem:elliptic-generator-criterion-e1-e2}
and~\ref{lem:o6-generator-criterion-iii}.}
\end{proof}

\subsection{\texorpdfstring{The case $n\ge3$}{The case n>=3}}

Throughout this subsection, assume $n\ge3$.
We prove that \textup{(E3)} does not hold, and hence
$E^-_{2n+2}(q)$ is not a triply even binary code.

\begin{fact}[{\cite{Ihringer}}]
\label{fac:ihringer-triple-common-neighbours-by-dimension}
Let $P_1,P_2,P_3\in\mathcal P_{2n+1}$ be distinct points satisfying
$P_1\sim P_2$, $P_2\sim P_3$, and $P_3\sim P_1$.
Then
{
\[
 \tau(P_1,P_2,P_3)=
 \begin{cases}
  \lambda-1,
    &\dim_{\mathbb F_q}(P_1+P_2+P_3)=2,\\[2mm]
  \begin{aligned}
   &q^2+q-2\\[-1mm]
   &\quad+q^3(q^{n-2}+1)\dfrac{q^{n-3}-1}{q-1},
  \end{aligned}
    &\dim_{\mathbb F_q}(P_1+P_2+P_3)=3.
 \end{cases}
\]
}
For details, see
\cite[Section~2.2, Lemma~4.1, and its proof]{Ihringer}.
\end{fact}

\begin{corollary}\label{cor:ihringer-triple-common-neighbours-parity}
Under the assumptions of
Fact~\ref{fac:ihringer-triple-common-neighbours-by-dimension}, we have
\begin{equation}\label{eq:tauis23}
 \tau(P_1,P_2,P_3)\equiv
 \begin{cases}
  1\pmod2,
    &\dim_{\mathbb F_q}(P_1+P_2+P_3)=2,\\
  0\pmod2,
    &\dim_{\mathbb F_q}(P_1+P_2+P_3)=3.
 \end{cases}
\end{equation}
\end{corollary}

\begin{proof}
By Fact~\ref{fac:elliptic-srg-parameters}, $\lambda$ is even.
Since $q$ is odd, the assertion holds.
\end{proof}

\begin{proposition}\label{lem:higher-rank-four-points}
There exist distinct
$P_1,P_2,P_3\in\mathcal P_{2n+1}\setminus\{P_0\}$ satisfying the following
conditions.
\begin{itemize}
 \item[\textup{(H1)}]
 $P_i\sim P_j$ for all $0\le i<j\le3$.

 \item[\textup{(H2)}]
 $\dim_{\mathbb F_q}(P_0+P_1+P_2)=2$, and
 $\dim_{\mathbb F_q}(P_0+P_1+P_3)
  =\dim_{\mathbb F_q}(P_0+P_2+P_3)
  =\dim_{\mathbb F_q}(P_1+P_2+P_3)=3$.
\end{itemize}
\end{proposition}

\begin{proof}
Let $M\subsetneq V$ be a maximal singular subspace containing $P_0$.
By the definition of an elliptic quadratic space,
$\dim_{\mathbb F_q}M=n\ge3$.
Hence we may choose a singular $2$-space $U$ and a singular $3$-space $T$
satisfying $P_0\subsetneq U\subsetneq T\subseteq M$.
Since $\ell_U$ is the set of $1$-dimensional subspaces of the $2$-space $U$,
it contains $q+1\ge4$ points. We may choose distinct
$P_1,P_2\in\ell_U\setminus\{P_0\}$, and then choose a point
$P_3\subsetneq T$ with $P_3\not\subseteq U$.

Since $T$ is singular, $P_i+P_j$ is a singular $2$-space for every
$0\le i<j\le3$. Hence $P_i$ and $P_j$ lie on the line
$\ell_{P_i+P_j}$, and therefore $P_i\sim P_j$.
Thus \textup{(H1)} holds.
By the choice of $P_0,P_1,P_2,$ and $P_3$, \textup{(H2)} also holds.
\end{proof}

\begin{corollary}\label{cor:higher-rank-not-triply-even}
$E^-_{2n+2}(q)$ is not a triply even binary code.
\end{corollary}

\begin{proof}
Take $P_1,P_2,P_3$ satisfying \textup{(H1)} and \textup{(H2)}.
{By Lemma~\ref{lem:triple-intersection-transfer} and
\eqref{eq:tauis23}},
{
\[
\begin{aligned}
 &\wt(g_{P_1}*g_{P_2}*g_{P_3})\\
 &\equiv \tau(P_1,P_2,P_3)+\tau(P_0,P_2,P_3)+\tau(P_0,P_1,P_3)+\tau(P_0,P_1,P_2)\\
 &\equiv1\pmod2.
\end{aligned}
\]
}
Thus \textup{(E3)} does not hold, and
Corollary~\ref{cor:elliptic-generator-criterion} shows that
$E^-_{2n+2}(q)$ is not a triply even binary code.
\end{proof}

\section{\texorpdfstring{{Maximality of $E^-_6(q)$}}{Maximality of E^-_6(q)}}

We prove that $E^-_6(q)$ is a maximal triply even binary code by showing
that \eqref{eq:elliptic-maximality-condition} holds for $n=2$ and $q\equiv3\pmod4$. For each $\ell\in\mathcal L_5$, let
$\one_\ell\in\mathbb F_2^{\mathcal P_5}$ denote its incidence vector, and put
$
 L^-_6(q)
:=\langle\one_\ell\mid\ell\in\mathcal L_5\rangle_{\mathbb F_2}$.

\begin{fact}[{\cite[Theorem~2.3]{BrouwerWilbrinkHaemers}}]
\label{fac:bwh-line-incidence-duality}
Let $q$ be an odd prime power.
Then we have
$
 E^-_6(q)^\perp=L^-_6(q)$,
 $\dim E^-_6(q)=q^3-q^2+q$, and $\dim L^-_6(q)=q^4+q^2+1$.
\end{fact}

For $P\in\mathcal P_5$, put
$\mathcal L(P):=\{\ell\in\mathcal L_5\mid P\in\ell\}$.
By Lemma~\ref{lem:o6-generalized-quadrangle},
$|\mathcal L(P)|=q^2+1$.

\begin{lemma}\label{lem:o6-line-sum}
Let $q$ be an odd prime power.
For every $P\in\mathcal P_5$, we have
$
 \sum_{\ell\in\mathcal L(P)}\one_\ell=r_P.
$
\end{lemma}

\begin{proof}
By (P1), every $P'\in N(P)$ lies on exactly one
line in $\mathcal L(P)$.
No point outside $N(P)\cup\{P\}$ lies on a line
in $\mathcal L(P)$.
Since $P$ lies on every line in $\mathcal L(P)$
and $q$ is odd, we have
\[
\left(\sum_{\ell\in\mathcal L(P)}\one_\ell\right)(P)
\equiv |\mathcal L(P)|
=(q^2+1)\equiv 0\pmod 2.
\]
Thus the claim follows.
\end{proof}

\begin{proposition}\label{prop:o6-line-schur-identity}
Let $q$ be an odd prime power, $\ell\in\mathcal L_5$ and $P\in\ell$.
For any distinct $Q_1,Q_2\in\ell\setminus\{P\}$,
we have
$(r_P+r_{Q_1})*(r_P+r_{Q_2})=r_P+\one_\ell$.
In particular, $r_P+\one_\ell\in E^-_6(q)*E^-_6(q)$
for every $\ell\in\mathcal L_5$ and every $P\in\ell$.
\end{proposition}

\begin{proof}
By \textup{(P2)$'$}, we have
$|\ell\setminus\{P\}|=q\ge3$.
Thus we may choose distinct points
$Q_1,Q_2\in\ell\setminus\{P\}$.

For $R\in\mathcal P_5$, let $\mathbf e_R$ be the vector whose
$R$-coordinate is $1$ and whose other coordinates are $0$.
By Lemma~\ref{lem:gq-triangles}, any two distinct points
$X,Y\in\ell$ satisfy
$N(X)\cap N(Y)=\ell\setminus\{X,Y\}$.
Thus $r_X*r_Y=\one_\ell+\mathbf e_X+\mathbf e_Y$.
Consequently,
\[
\begin{aligned}
 &r_P*r_{Q_1}+r_P*r_{Q_2}+r_{Q_1}*r_{Q_2}\\
 &=(\one_\ell+\mathbf e_P+\mathbf e_{Q_1})
  +(\one_\ell+\mathbf e_P+\mathbf e_{Q_2})+(\one_\ell+\mathbf e_{Q_1}+\mathbf e_{Q_2})=\one_\ell.
\end{aligned}
\]
Since $r_P*r_P=r_P$, we obtain
\[
 (r_P+r_{Q_1})*(r_P+r_{Q_2})
 =r_P+r_P*r_{Q_2}+r_P*r_{Q_1}
   +r_{Q_1}*r_{Q_2}
 =r_P+\one_\ell.
\]
Since $r_P+r_{Q_i}\in E^-_6(q)$ for $i=1,2$,
the last assertion follows.
\end{proof}

\begin{theorem}\label{thm:o6-square-dual-line}
Let $q$ be an odd prime power with $q\equiv3\pmod{4}$. Then we have
$
 E^-_6(q)*E^-_6(q)=E^-_6(q)^\perp=L^-_6(q)$, and
 $\dim\bigl(E^-_6(q)*E^-_6(q)\bigr)=q^4+q^2+1$.
\end{theorem}

\begin{proof}
Let $P\in\mathcal P_5$.
By Proposition~\ref{prop:o6-line-schur-identity},
we have
$r_P+\one_\ell\in E^-_6(q)*E^-_6(q)$
for every $\ell\in\mathcal L(P)$.
Since $|\mathcal L(P)|=q^2+1$ is even,
Lemma~\ref{lem:o6-line-sum} gives
\[
 \sum_{\ell\in\mathcal L(P)}(r_P+\one_\ell)
 =(q^2+1)r_P+r_{P}=r_P.
\]
Hence $r_P\in E^-_6(q)*E^-_6(q)$. 

Now let $\ell\in\mathcal L_5$ and $P\in\ell$.
By Proposition~\ref{prop:o6-line-schur-identity},
\[
 \one_\ell=(r_P+\one_\ell)+r_P
 \in E^-_6(q)*E^-_6(q).
\]
Thus $L^-_6(q)\subseteq E^-_6(q)*E^-_6(q)$.

Proposition~\ref{prop:triply-even-square-maximal}\textup{(M1)},
Proposition \ref{prop:o6-triply-even}, and
Fact~\ref{fac:bwh-line-incidence-duality} give
$
 E^-_6(q)*E^-_6(q)=E^-_6(q)^\perp=L^-_6(q)$ and $\dim\bigl(E^-_6(q)*E^-_6(q)\bigr)=q^4+q^2+1$.
\end{proof}

\begin{corollary}\label{cor:o6-maximal}
Let $q$ be an odd prime power with $q\equiv3\pmod{4}$. Then $E^-_6(q)$ is a
$
 \bigl[(q+1)(q^3+1),\ q(q^2-q+1)\bigr]
$
maximal triply even binary code.
\end{corollary}

\begin{proof}
By Corollary~\ref{cor:elliptic-maximality-criterion}, Proposition~\ref{prop:o6-triply-even}, and
Theorem~\ref{thm:o6-square-dual-line},
$E^-_6(q)$ is a maximal triply even binary code.
Moreover, \eqref{eq:o6-srg-parameters} and
Fact~\ref{fac:bwh-line-incidence-duality} give
$|\mathcal P_5|=(q+1)(q^3+1)$ and
$\dim E^-_6(q)=q(q^2-q+1)$, respectively.
\end{proof}

\section*{Acknowledgments}
This work was supported by JST SPRING, Grant Number JPMJSP2152.
We thank Professor Koichi Betsumiya for his guidance and encouragement
throughout this research.

\end{document}